\documentclass[a4paper,11pt]{amsart}
\usepackage{amsmath, amsfonts, amssymb, amsthm, amscd}
\usepackage[colorlinks=true, linkcolor=blue, citecolor=blue]{hyperref}
\usepackage{graphicx}
\usepackage{enumerate}
\usepackage{url}
\usepackage{color}
\usepackage[utf8]{inputenc}
\usepackage{tikz}
\usetikzlibrary{arrows}
\usetikzlibrary{patterns} 
\usepackage[a4paper,scale={0.72,0.74},marginratio={1:1},footskip=7mm,headsep=10mm]{geometry}

\numberwithin{equation}{section}

\newtheorem{theorem}{Theorem}

\newtheorem{proposition}[theorem]{Proposition}

\newtheorem{remark}[theorem]{Remark}

\newcommand{\bbE}{{\ensuremath{\mathbb E}} }

\newcommand{\bbN}{{\ensuremath{\mathbb N}} }

\newcommand{\bbP}{{\ensuremath{\mathbb P}} }

\newcommand{\bbR}{{\ensuremath{\mathbb R}} }

\newcommand{\bbZ}{{\ensuremath{\mathbb Z}} }

\newcommand{\cE}{{\ensuremath{\mathcal E}} }

\newcommand{\cP}{{\ensuremath{\mathcal P}} }

\newcommand{\cZ}{{\ensuremath{\mathcal Z}} }

\newcommand{\ga}{\alpha}
\newcommand{\gb}{\beta}

\newcommand{\gD}{\Delta}
\newcommand{\gep}{\varepsilon}

\newcommand{\gl}{\lambda}

\newcommand{\go}{\omega}

\renewcommand{\tilde}{\widetilde}          
\DeclareMathSymbol{\leqslant}{\mathalpha}{AMSa}{"36} 
\DeclareMathSymbol{\geqslant}{\mathalpha}{AMSa}{"3E} 
\DeclareMathSymbol{\eset}{\mathalpha}{AMSb}{"3F}     

\newcommand\bP{\ensuremath{\mathrm{P}}}
\newcommand\bE{\ensuremath{\mathrm{E}}}

\renewcommand{\epsilon}{\varepsilon}

\newcommand{\ind}{{\sf 1}}

\newcommand{\card}{\mathrm{card}}

\newenvironment{myenumerate}{
\renewcommand{\theenumi}{\arabic{enumi}}
\renewcommand{\labelenumi}{{\rm(\theenumi)}}
\begin{list}{\labelenumi}
{
\setlength{\itemsep}{0.4em}
\setlength{\topsep}{0.5em}
\setlength\leftmargin{2.45em}
\setlength\labelwidth{2.05em}
\setlength{\labelsep}{0.4em}
\usecounter{enumi}
}
}
{\end{list}
}

\renewenvironment{enumerate}{
\begin{myenumerate}}
{\end{myenumerate}}

\newcommand{\beq}{\begin{equation}}
\newcommand{\eeq}{\end{equation}}
\newcommand{\ba}{\begin{aligned}}
\newcommand{\ea}{\end{aligned}}

\begin{document}

\title[]{A connection between the random pinning model and random walks in sparse random environments}
\author{Julien Poisat}
\address{Université Paris-Dauphine, CNRS, UMR 7534, CEREMADE, PSL Research University, Place du Maréchal de Lattre de Tassigny, 75016 Paris, France.}
\email{poisat@ceremade.dauphine.fr}

\begin{abstract}
The purpose of this short note is to establish a connection between a one-dimensional random walk in a random sparse environment and the random pinning model. We show that the grand canonical partition function of the pinning model coincides with the mean number of returns to the origin for a random walk in a random sparse environment averaged over the randomness location. We obtain thereof some information on the integrability of the number of return times in the annealed and partially annealed setups.
\end{abstract}

\thanks{The author acknowledges the support of French grant ANR LOCAL (ANR-22-CE40-0012).}
 
\date{\today}

\maketitle

\tableofcontents

\section{Introduction}

Random walk in a \emph{sparse} random environment (RWsRE) is a variant of the more common random walk in random environment (RWRE) that behaves as a simple random walk everywhere but at \emph{random locations} where random transition probabilities are drawn. This model with a double layer of randomness has recently received some interest~\cite{BDIM20,BDIMR19,BDK24,Dag17,MRS16,Seol13}. In this short paper we establish a connection between a one-dimensional RWsRE and the so-called random pinning model~\cite{Gia07,Gia11}. We first remind the reader of standard facts about the random walk in a general potential and in a random environment, in Sections~\ref{sec:rw-pot} and~\ref{sec:rwre} respectively. We then recall some (non-exhaustive) results on the RWsRE in Section~\ref{sec:rwsre} and finally make the connection with the random pinning model in Section~\ref{sec:pinning}. The key relation, displayed in Proposition~\ref{pr:key-identity}, is an equality between the grand canonical partition function of the random pinning model and the mean number of return times to the origin for the law of the RWsRE averaged over the \emph{location} of the randomness only (partially annealed setup). Our observations eventually lead to Theorem~\ref{thm:rwrse-thresholds}: the annealed and (more interestingly) quenched critical curves of the pinning model are shown to correspond to thresholds for the (respectively annealed and partially annealed) integrability of the mean number of return times of the RWsRE.\\

{\bf \noindent Notation.} Throughout the paper, $\bbN_0$ and $\bbN$ respectively denote the set of non-negative and positive integers.

\section{Reminder on random walks in a potential}
\label{sec:rw-pot}
Let us first recall a few standard facts about random walks in a potential. Let $(V_i)_{i\in \bbZ}$ be a real-valued sequence with $V_0=0$ and $\gD V_i := V_i - V_{i-1}$ for $i\in \bbZ$. The \emph{random walk in potential $V$} is the Markov chain $X= (X_n)_{n\ge 0}$ on $\bbZ$, started at $X_0=0$ unless stated otherwise, with transition matrix:
\beq
\label{eq:trans-prob}
p(i,i+1)= 1-p(i,i-1) = \frac{1}{1+e^{\gD V_i}}\in (0,1), \qquad i\in \bbZ.
\eeq
Its law shall be denoted by $\cP_V$ in the sequel, with $\cE_V$ the corresponding expectation.
Alternatively to~\eqref{eq:trans-prob},
\beq
\gD V_i = \log\Big[ \frac{p(i,i-1)}{p(i,i+1)} \Big], \qquad i\in \bbZ,
\eeq
and one may also think of $e^{-V_i}$ as a conductance assigned to the edge $(i,i+1)$.
Set $W(0) = 0$ and
\beq
W(n) := \sum_{0\le k<n} e^{V_k},\qquad
W(-n) := -\sum_{1\le k\le n} e^{V_{-k}}
 \qquad n\ge1.
\eeq
The following fact is standard: using that $(W(X_n))_{n\ge 0}$ is a martingale adapted to $X$, an application of Doob's optional stopping theorem yields
\beq
\cP_V( H_a > \tilde H_c| X_0 = b) = \frac{W(b)-W(a)}{W(c)-W(a)}, \qquad a< b\le c,
\eeq
where
\beq
\ba
\tilde H_x &:= \inf\{n\ge 0\colon X_n = x\}\quad {\rm (hitting\ time),}\\
H_x &:= \inf\{n\ge 1\colon X_n = x\}\quad {\rm (entrance\ time),}
\ea
\eeq
see e.g.~\cite{Zei2012}.
In particular,
\beq
\cP_V(H_y > \tilde H_x| X_0 = y+1)^{-1} = \sum_{y\le i< x} e^{V_i - V_y}, \qquad 0\le y < x.
\eeq
Choosing $y=0$ and letting $x\to\infty$ in the line above, we obtain
\beq
\label{eq:escfromone}
\cP_V(H_0 = +\infty| X_0 = 1)^{-1} = \sum_{i\ge 0} e^{V_i},
\eeq
and in a similar way,
\beq
\label{eq:escfromminusone}
\cP_V(H_0 = +\infty| X_0 = -1)^{-1} = \sum_{i\ge 0} e^{V_{-i-1} - V_{-1}}.
\eeq
If we assume moreover that $V$ is symmetric w.r.t to $-1/2$, so that $V_i = V_{-i-1} - V_{-1}$ for every $i\ge 0$, then the quantities in~\eqref{eq:escfromone} and~\eqref{eq:escfromminusone} coincide. Consequently, 
we obtain the following expression for the mean number of return times to the origin:
\beq
\label{eq:return-times}
\cE_V( \card\{n\ge 0 \colon X_n = 0\}) = \sum_{i\ge 0} e^{V_i}.
\eeq
In particular, $X$ is \emph{transient} (resp. \emph{recurrent}) if and only if the above sum converges (resp. diverges). Throughout the paper, we assume the aforementionned symmetry on the potential $V$, which amounts to replacing the $\bbZ$-valued Markov chain by an $\bbN_0$-valued chain (a.k.a.\ birth-and-death process).
\section{Random walk in a random environment (RWRE)}
\label{sec:rwre}
A random walk in a random environment is obtained when the potential is drawn {\it randomly}. This is a very active field of research, see~\cite{DevDiel2023,DreRam14,Hug95,Hug96,Kesten75,Rev75,Sinai82,Solomon75,Sznitman04,Zei2012} for a selection of seminal papers and reviews on the topic. A common assumption made in the literature is that the increments $(\gD V_i)_{i\ge 1}$ form an i.i.d. sequence. Therefore, we let $(\go_i)_{i\ge 1}$ be a sequence of (integrable) i.i.d. random variables with law $\bbP$ and define
\beq
\gD V_i := \gb \go_i + h, \qquad i\ge 1,
\eeq 
where $h\in \bbR$ and (w.l.o.g.) $\gb\ge 0$ and $\bbE(\go_1) = 0$. In that setup, it is known since Solomon~\cite{Solomon75} that $X$ is transient under $\cP_V$ if and only if $h<0$\footnote{The reader might be puzzled here: Solomon's criterion states that, provided that the expectation of $\gD V_i$ is defined, transience holds if and only if $\bbE(\gD V_i)\neq 0$, that is \emph{$h\neq 0$}. However, due to the symmetry assumption made on the potential $V$ at the end of Section~\ref{sec:rw-pot}, the condition for transience becomes $h<0$ in our case.} and that $X$ is transient with \emph{positive} speed if and only if $h< -\log \bbE(e^{\gb \go_1})$ (the exponential moment may be finite or infinite). Remark that, by Jensen's inequality, $-\log \bbE(e^{\gb \go_1})\le 0$ with the inequality being strict if $\go_1$ is not identically zero.
\section{Random walk in a sparse random environment (RWsRE)}
\label{sec:rwsre}
More recently, several authors~\cite{BDIM20,BDIMR19,BDK24,Dag17,MRS16,Seol13} have introduced a variant of RWRE, by assuming that the randomness of the potential is only seen at \emph{random locations}, outside which the potential is flat. This means defining the potential by $V_0 = 0$ and
\beq
\label{eq:rwsre-pot}
\gD V_i = (h+\gb\go_{i}) \ind_{\{i\in \tau\}}, \qquad i\ge 1,
\eeq
where $\tau = \{\tau_i\}_{i\ge 0}$ is a certain random subset of $\bbN_0$. A common assumption, that we shall adopt in the sequel, is that $(\tau_i)_{i\ge 0}$ is a \emph{renewal} process starting at $\tau_0=0$ with independent $\bbN$-valued increments $(\tau_{i}-\tau_{i-1})_{i\ge 1}$. Its law shall be denoted by $\bP$. In that case, we write
\beq
\cP_{\gb,h}^{\go, \tau}:= \text{law of the quenched RWsRE when \eqref{eq:rwsre-pot} holds,}
\eeq
instead of $\cP_V$. Also, we assume here that $\tau$ and $\go$ are independent, which is not necessarily the case in the RWsRE literature. Such a random environment is called \emph{moderately} (resp. \emph{strongly}) sparse if $\bE(\tau_1)$ is finite (resp. infinite). By a result of Matzavinos, Roitershtein and Seol~\cite[Theorem 3.1]{MRS16}, we have the following:
\begin{proposition}
\label{pr:ann-transient}
Assume that $\bE(\log \tau_1) < \infty$. Then the RWsRE is transient under the \emph{annealed} law, i.e.\ 
\beq
\label{eq:ann-transient}
\bbE\bE\Big[\cP_{\gb,h}^{\go, \tau}(\card\{n\ge 0 \colon X_n = 0\} < \infty)\Big] = 1,
\eeq
whenever $h<0$\footnote{Same remark as above.}.
\end{proposition}
In complement to Proposition~\ref{pr:ann-transient}, the same authors~\cite[Theorem 3.3]{MRS16} proved under the same assumptions, and still considering the annealed law, that the RWsRE is transient with \emph{positive} speed if $h < -\log \bbE(e^{\gb \go_1})$ and $\bE(\tau_1^2) < \infty$, or transient with \emph{zero} speed otherwise. See~\cite[Proposition 2.1]{BDIMR19} in the case where the vectors $(\go_k, \tau_k - \tau_{k-1})_{k\ge 1}$ are independent and identically distributed (i.i.d.) but $\go_k$ and $\tau_k - \tau_{k-1}$ are allowed to be dependent.

\par Note that the equality in~\eqref{eq:ann-transient} readily implies that
\beq
\cP_{\gb,h}^{\go, \tau}(\card\{n\ge 0 \colon X_n = 0\} < \infty) = 1, \qquad \bP\otimes\bbP-{\rm a.s.\ ,}
\eeq
that is transience for almost every quenched sparse environment. This also implies that
\beq
\cE_{\gb,h}^{\go, \tau}(\card\{n\ge 0 \colon X_n = 0\}) < \infty, \qquad \bP\otimes\bbP-{\rm a.s.}
\eeq
However,~\eqref{eq:ann-transient} does \emph{not} imply that 
\beq
\bbE\bE\cE_{\gb,h}^{\go, \tau}(\card\{n\ge 0 \colon X_n = 0\}) < \infty.
\eeq
We will see in Section~\ref{sec:pinning} that the mean number of return times of RWsRE to the origin, once averaged over the renewal times, coincides with the grand canonical partition function of the so-called \emph{random pinning model}. This will allow us to investigate the integrability of the former quantity w.r.t. the law $\bP\otimes \bbP$ and (more interestingly) the law $\bP$, for $\bbP$-a.e.\ sequence $\go$. We will finally show in Theorem~\ref{thm:rwrse-thresholds}, that is our main result, that there exist several regimes of integrability depending on the value of $h$, and that these regimes happen to be delimited by the (quenched and annealed) critical points of the random pinning model.
\section{Connection with the random pinning model}
\label{sec:pinning}
Let us now introduce the random pinning model, to which we want to connect the RWsRE. This is a statistical mechanics model of a polymer pinned along a (one-dimensional) defect line that is motivated for instance by interface models or the phenomenon of DNA denaturation (a.k.a.\ Poland-Scheraga models in this context). There has been an intense and fruitful activity on this model lately, in particular around the issue of disorder relevance, which we shortly come back to below. We keep our exposition as concise as possible but invite the reader to refer to~\cite{Gia07,Gia11} or \cite[Section 11]{dH09},  and references therein for a more thorough introduction to this model. From now on we will assume that
\beq
\bP(\tau_1 = n) = L(n)n^{-(1+\ga)}, \qquad n\in \bbN,
\eeq
where $L$ is a slowly varying function and $\ga\ge 0$, as well as
\beq
\bbE(\go_1) = 0, \qquad \bbE(e^{\gb \go_1}) < \infty, \qquad \forall\gb\ge0.
\eeq
The quenched \emph{partition function} of the random pinning model is defined as
\beq
Z_{n,\go}^{\gb,h} := \bE\Big[\exp\Big(\sum_{i=1}^n (\gb \go_i + h)\ind_{\{i\in \tau\}}\Big) \Big], \qquad n\ge 1,
\label{eq:pf}
\eeq
with the convention $Z_{0,\go}^{\gb,h} := 0$. Here, the renewal set $\tau$ can be thought of as the random set of contact points between a linear polymer chain and a one-dimensional defect line, with $\gb \go_i + h$ the random energy reward/penalty gathered from a contact at site $i\in \bbN$. It is known that the following limit, called quenched \emph{free energy}, exists $\bbP$-a.s.\ 
\beq
f(\gb,h) := \lim_{n\to\infty} \frac 1n \log Z_{n,\go}^{\gb,h} \ge 0,
\eeq
and that the following dichotomy holds: there exists a critical curve $-\infty <h_c(\gb)\le 0$ such that
\beq
\label{eq:fe-dichotomy}
f(\gb, h) = 0\,\,\text{for every $h \le h_c(\gb)$}\qquad \text{and} \qquad f(\gb, h) >0\,\,\text{for every $h>h_c(\gb)$}.
\eeq
This is a signature of a phase transition between a delocalized phase and a localized phase, and it is known that that (i) $h_c(0)= 0$ and (ii) $h_c(\gb)<0$ if $\gb>0$, see~\cite{AS06} or \cite[Chap. 5, Sect. 2]{Gia07}. One of the major challenges in this model has been to determine whether disorder is \emph{relevant} or not, that is whether the critical curve coincides with the critical curve for the \emph{annealed} model for at least sufficiently small values of $\gb$ (disorder strength). Note that the annealed critical curve, denoted by $h_c^{\rm a}(\gb)$, is easy to compute, since the annealed partition function equals
\beq
\label{eq:ann-model}
\bbE Z_{n,\go}^{\gb,h} = \bE\Big[\exp\Big(\sum_{i=1}^n (\gl(\gb) + h)\ind_{\{i\in \tau\}}\Big) \Big], \qquad \text{where}\quad \gl(\gb) := \log \bbE(e^{\gb\go_1}).
\eeq
Therefore,
\beq
\label{eq:ann-model2}
h_c^{\rm a}(\gb) = -\gl(\gb)
\eeq
is simply the critical point of an adequately shifted homogeneous (i.e.\ $\gb=0$) model. The following equivalence for disorder relevance in the sense of critical point shifts is now known to hold: 
\beq
\{\forall \gb>0,\, h_c(\gb) > h_c^{\rm a}(\gb)\}\quad
\Leftrightarrow \quad
\sum_{n\ge 1} \frac{1}{n^{2(1-\ga)}L(n)^2} = +\infty,
\label{eq:relevance-crit}
\eeq
that is, if and only if the renewal process $\tau^{(1)}\cap \tau^{(2)}$, obtained by intersecting two independent copies of $\tau$, is transient, see~\cite[Theorem 2.2]{BerLac18} and references therein. The result displayed in \eqref{eq:relevance-crit} finally put on firm mathematical grounds a prediction made by physicists (the Harris criterion) and resolved a debate about the delicate marginal case $\ga = 1/2$. Let us stress however that the strict inequality between the quenched and annealed critical points always holds for large enough $\gb$~\cite[Corollary 3.2]{Ton-AAP2008} unless $\ga=0$, in which case $h_c(\gb) = h_c^{\rm a}(\gb)$ for every $\gb\ge 0$~\cite{AlexZyg2010}.

\par Before getting to our main observation, we define the \emph{grand canonical} partition function of the model by
\beq\label{eq:gcpf}
\cZ_{\gb,h}^\go(f) := \sum_{n\ge 0} Z_{n,\go}^{\gb,h} e^{-fn}, \qquad f\in \bbR.
\eeq
Note that
\beq\label{eq:fe2}
f(\gb,h) = \inf\{f\ge 0 \colon \cZ_{\gb,h}^\go(f)< \infty\ \bbP{\rm -a.s.}\}.
\eeq
It turns out that this is related to the RWsRE with an additional drift. To this end, we let
\beq
\label{eq:def-cP-gb-h-f}
\cP_{\gb,h,f}^{\go,\tau} := \cP_V \quad {\rm when} \quad 
\gD V_i := (h+\gb\go_{i})\ind_{\{i\in \tau\}}-f, \quad i\ge 1.
\eeq
The following identity, which connects the RWsRE to the random pinning model, is the key for the remaining discussion:
\begin{proposition}
\label{pr:key-identity}
For every $\gb\ge 0$ and $h,f\in\bbR$,
\beq
\label{eq:key_rel}
\bE\cE_{\gb,h,f}^{\go,\tau}( \card\{n\ge 0 \colon X_n = 0\}) = \cZ_{\gb,h}^\go(f).
\eeq
\end{proposition}
\begin{proof}[Proof of Proposition~\ref{pr:key-identity}]
Using~\eqref{eq:return-times} in combination with~\eqref{eq:def-cP-gb-h-f}, one obtains
\beq
\cE_{\gb,h,f}^{\go,\tau}( \card\{n\ge 0 \colon X_n = 0\}) = \sum_{n\ge 0} \exp\Big(\sum_{1\le i \le n} (h+\gb\go_{i})\ind_{\{i\in \tau\}} - fn\Big),
\eeq
where the sum in the exponential is understood as $0$ when $n=0$. The claim now follows from~\eqref{eq:pf} and~\eqref{eq:gcpf}, once we average over $\tau$ the above identity.
\end{proof}
From~\eqref{eq:fe-dichotomy}, \eqref{eq:fe2} and Proposition~\ref{pr:key-identity} we readily obtain the following:
\begin{proposition}
\label{pr:dichotomy}
Let $\gb \ge0$.
\begin{enumerate}
\item If $h> h_c(\gb)$ then $\cZ_{\gb,h}^\go(\gep)=\bE\cE_{\gb,h,\gep}^{\go,\tau}( \card\{n\ge 0 \colon X_n = 0\})=\infty$ for every $0\le \gep < f(\gb,h)$.
\item If $h\le h_c(\gb)$ then $\cZ_{\gb,h}^\go(\gep)= \bE\cE_{\gb,h,\gep}^{\go,\tau}( \card\{n\ge 0 \colon X_n = 0\})< \infty$ for every $\gep>0$.
\end{enumerate}
\end{proposition}
It is however not straightforward to deduce a result for $\cP^{\go,\tau}_{\gb,h,f}$ in the  case $f=0$ and $h\le h_c(\gb)$ from the mere definition of the free energy. Fortunately, we may use a result from Mourrat~\cite[Theorem 1]{Mourrat12} about the \emph{pinned} (or \emph{constrained}) version of the grand canonical partition function, that is
\beq
\cZ_{\gb,h}^{\go,{\rm c}}(f) := \sum_{n\ge 0} Z_{n,\go}^{\gb,h}(n\in \tau) e^{-fn}, \qquad f\in \bbR,
\eeq
where
\beq
Z_{n,\go}^{\gb,h}(n\in \tau) := \bE\Big[\exp\Big(\sum_{i=1}^n (\gb \go_i + h)\ind_{\{i\in \tau\}}\Big) \ind_{\{n\in \tau\}} \Big], \quad n\ge 1, \qquad (0 \text{ if } n=0),
\eeq
is the \emph{pinned} (or \emph{constrained}) partition function (compare to~\eqref{eq:pf}). With this notation in hand, 
\beq
\label{eq:mourrat}
\forall \gb\ge 0,\qquad
h< h_c(\gb)\quad
\Rightarrow \quad
\cZ_{\gb,h}^{\go,{\rm c}}(0) < \infty\quad \text{$\bbP$-a.s.},
\eeq
according to~\cite[Theorem 1]{Mourrat12} (mind the minus sign in front of $h$ therein). We can now relate the grand canonical partition function to its pinned version. By decomposing the partition function according to the last renewal point before $n$, we obtain
\beq
Z_{n,\go}^{\gb,h} = \sum_{k\ge 0} \bE\Big(e^{\sum_{1\le i \le k} (\gb \go_{\tau_i} + h)} 
\ind_{\{\tau_k \le n < \tau_{k+1}\}}\Big),
\eeq
(the sum in the above exponential is understood as zero if $k=0$).
Therefore,
\beq
\label{eq:p-gcpf-to-f-gcpf}
\ba
\cZ_{\gb,h}^{\go}(0) = \sum_{n\ge 0} Z_{n,\go}^{\gb,h} &= 
\bE\Big[ \sum_{k\ge 0} e^{\sum_{1\le i \le k} (\gb \go_{\tau_i} + h)} (\tau_{k+1}-\tau_k) \Big]\\
 &=  \bE(\tau_1)
\sum_{k\ge 0} \bE\Big(e^{\sum_{1\le i \le k} (\gb \go_{\tau_i} + h)}\Big)
= \bE(\tau_1) \cZ_{\gb,h}^{\go,{\rm c}}(0).
\ea
\eeq
We deduce thereby:
\begin{proposition}
\label{pr:dichotomy2}
The following dichotomy holds:
\begin{enumerate}
\item If $\bE(\tau_1) = \infty$ then $\cZ_{\gb,h}^{\go}(0)= \infty$ $(\bbP\text{-a.s.})$ for every $\gb \ge 0$ and $h\in \bbR$.
\item If $\bE(\tau_1) < \infty$ then $\cZ_{\gb,h}^{\go}(0)< \infty$ $(\bbP\text{-a.s.})$ for every $\gb\ge0$ and $h< h_c(\gb)$.
\end{enumerate}
\end{proposition}
\begin{proof}[Proof of Proposition~\ref{pr:dichotomy2}]
Item (1) follows from~\eqref{eq:p-gcpf-to-f-gcpf} while Item (2) follows from~\eqref{eq:p-gcpf-to-f-gcpf} and~\eqref{eq:mourrat}.
\end{proof}
\begin{remark}
Note that the first item in the above dichotomy may be obtained in a different manner. Indeed, by a standard ruin probability estimate,
\beq
\cP_{\gb,h}^{\go,\tau}(H_0 = \infty) \le \cP_{\gb,h}^{\go,\tau}(H_0 > H_{\tau_1}) = 1/\tau_1,
\eeq
which leads, using Proposition~\ref{pr:key-identity}, to
\beq
\cZ_{\gb,h}^\go(0) = \bE\cE_{\gb,h}^{\go,\tau}( \card\{n\ge 0 \colon X_n = 0\}) \ge \bE(\tau_1).
\eeq
\end{remark}
Summarizing all what precedes, we obtain the three following cases:
\begin{theorem}
\label{thm:rwrse-thresholds}
Let $\gb>0$.
\begin{enumerate}
\item If $h_c(\gb)<h<0$ and $\bE(\log \tau_1) < \infty$ then the RWsRE is transient, so $\cE_{\gb,h}^{\go,\tau}( \card\{n\ge 0 \colon X_n = 0\})$ is finite $\bP\otimes \bbP$-a.s, but $\bE\cE_{\gb,h,\gep}^{\go,\tau}( \card\{n\ge 0 \colon X_n = 0\})$ is infinite $\bbP$-a.s.\ for every $0\le \gep < f(\gb,h)$.
\item If $h_c^{\rm a}(\gb) < h \le h_c(\gb)$ then $\bE\cE_{\gb,h,\gep}^{\go,\tau}( \card\{n\ge 0 \colon X_n = 0\})$ is finite $\bbP$-a.s.\ for every $\gep>0$. The result also holds for $\gep =0$ if $\bE(\tau_1)$ is finite and $h< h_c(\gb)$. However,
$\bbE\bE\cE_{\gb,h,\gep}^{\go,\tau}( \card\{n\ge 0 \colon X_n = 0\})$ is infinite for every $0\le \gep < f(0,h-h_c^{\rm a}(\gb))$.
\item If $h\le  h_c^{\rm a}(\gb)$ then $\bbE\bE\cE_{\gb,h,\gep}^{\go,\tau}( \card\{n\ge 0 \colon X_n = 0\})$ is finite for every $\gep>0$. The result also holds for $\gep =0$ if $\bE(\tau_1)$ is finite and $h< h_c^{\rm a}(\gb)$.
\end{enumerate}
\end{theorem}
Case (2) in Theorem~\ref{thm:rwrse-thresholds} is non-empty (i) for large enough $\gb$ when $\ga >0$ or (ii) for every $\gb>0$ if $\ga>1/2$ or $\ga = 1/2$ and $\sum n^{-1}L(n)^{-2} = +\infty$, see the discussion around~\eqref{eq:relevance-crit}.
\begin{proof}[Proof of Theorem~\ref{thm:rwrse-thresholds}]
\begin{enumerate}
\item The first part of the statement follows from Proposition~\ref{pr:ann-transient} while the second part follows from Item (1) in Proposition~\ref{pr:dichotomy}.
\item The first part of the statement follows from Item (2) in Proposition~\ref{pr:dichotomy}, while the second part follows from Item (2) in Proposition~\ref{pr:dichotomy2}. Averaging over $\go$ the identity in Proposition~\ref{pr:key-identity} and using~\eqref{eq:ann-model}-\eqref{eq:ann-model2}, we obtain
\beq
\label{eq:ann-identity}
\bbE\bE\cE_{\gb,h,\gep}^{\go,\tau}( \card\{n\ge 0 \colon X_n = 0\}) = \bbE \cZ_{\gb,h}^\go(f) = \cZ_{0,h-h_c^{\rm a}(\gb)}^\go(f).
\eeq
Combined with Item (1) in Proposition~\ref{pr:dichotomy} and the fact that $h_c(0)=0$, this proves the last part of the statement.
\item The first part of the statement follows from~\eqref{eq:ann-identity} combined with Item (2) in Proposition~\ref{pr:dichotomy} and the fact that $h_c(0)=0$. The second part of the statement follows from Item (2) in Proposition~\ref{pr:dichotomy2}.
\end{enumerate}
\end{proof}
\begin{remark}
The relation in Proposition~\ref{pr:key-identity} could be written for other one-dimensional statistical mechanics models based on renewal processes, like the copolymer model~\cite[Section 10]{dH09} or to the random pinning model with other types of disorder sequences, like~\emph{correlated} random sequences~\cite{Berger2014,CheChiPoi19}. In particular, \cite[Theorem 2.5]{Berger2013} asserts that when $\go$ is a Gaussian sequence with \emph{non-summable} non-negative correlations, then the free energy of the random pinning model is positive for every $\gb>0$ and $h\in\bbR$, i.e. $h_c(\gb) = -\infty$. Consequently, in that case $\bE\cE_{\gb,h}^{\go,\tau}( \card\{n\ge 0 \colon X_n = 0\})$ is infinite $\bbP$-a.s.\ for every $\gb>0$, \emph{no matter how small} $h$ is.
\end{remark}

\bibliographystyle{abbrv}
\bibliography{biblio.bib}

\begin{thebibliography}{10}

\bibitem{AS06}
K.~S. Alexander and V.~Sidoravicius.
\newblock Pinning of polymers and interfaces by random potentials.
\newblock {\em Ann. Appl. Probab.}, 16(2):636--669, 2006.

\bibitem{AlexZyg2010}
K.~S. Alexander and N.~Zygouras.
\newblock Equality of critical points for polymer depinning transitions with
  loop exponent one.
\newblock {\em Ann. Appl. Probab.}, 20(1):356--366, 2010.

\bibitem{Berger2013}
Q.~Berger.
\newblock Comments on the influence of disorder for pinning model in correlated
  {G}aussian environment.
\newblock {\em ALEA Lat. Am. J. Probab. Math. Stat.}, 10(2):953--977, 2013.

\bibitem{Berger2014}
Q.~Berger.
\newblock Pinning model in random correlated environment: appearance of an {\it
  infinite disorder} regime.
\newblock {\em J. Stat. Phys.}, 155(3):544--570, 2014.

\bibitem{BerLac18}
Q.~Berger and H.~Lacoin.
\newblock Pinning on a defect line: characterization of marginal disorder
  relevance and sharp asymptotics for the critical point shift.
\newblock {\em J. Inst. Math. Jussieu}, 17(2):305--346, 2018.

\bibitem{BDIM20}
D.~Buraczewski, P.~Dyszewski, A.~Iksanov, and A.~Marynych.
\newblock Random walks in a strongly sparse random environment.
\newblock {\em Stochastic Process. Appl.}, 130(7):3990--4027, 2020.

\bibitem{BDIMR19}
D.~Buraczewski, P.~Dyszewski, A.~Iksanov, A.~Marynych, and A.~Roitershtein.
\newblock Random walks in a moderately sparse random environment.
\newblock {\em Electron. J. Probab.}, 24:Paper No. 69, 44, 2019.

\bibitem{BDK24}
D.~Buraczewski, P.~Dyszewski, and A.~Ko\l~odziejska.
\newblock Weak quenched limit theorems for a random walk in a sparse random
  environment.
\newblock {\em Electron. J. Probab.}, 29:Paper No. 7, 30, 2024.

\bibitem{CheChiPoi19}
D.~Cheliotis, Y.~Chino, and J.~Poisat.
\newblock The random pinning model with correlated disorder given by a renewal
  set.
\newblock {\em Ann. H. Lebesgue}, 2:281--329, 2019.

\bibitem{Dag17}
K.~Dagtoros.
\newblock {\em Large {D}eviation {R}esults for {R}andom {W}alks in a {S}parse
  {R}andom {E}nvironment}.
\newblock ProQuest LLC, Ann Arbor, MI, 2017.
\newblock Thesis (Ph.D.)--Iowa State University.

\bibitem{dH09}
F.~den Hollander.
\newblock {\em Random polymers}, volume 1974 of {\em Lecture Notes in
  Mathematics}.
\newblock Springer-Verlag, Berlin, 2009.
\newblock Lectures from the 37th Probability Summer School held in Saint-Flour,
  2007.

\bibitem{DevDiel2023}
A.~Devulder, R.~Diel, and X.~Zeng.
\newblock Some recent advances in random walks and random environments.
\newblock In {\em Journ\'ees {MAS} 2020---random modelization and physics},
  volume~74 of {\em ESAIM Proc. Surveys}, pages 38--50. EDP Sci., Les Ulis,
  2023.

\bibitem{DreRam14}
A.~Drewitz and A.~F. Ram\'irez.
\newblock Selected topics in random walks in random environment.
\newblock In {\em Topics in percolative and disordered systems}, volume~69 of
  {\em Springer Proc. Math. Stat.}, pages 23--83. Springer, New York, 2014.

\bibitem{Gia07}
G.~Giacomin.
\newblock {\em Random polymer models}.
\newblock Imperial College Press, London, 2007.

\bibitem{Gia11}
G.~Giacomin.
\newblock {\em Disorder and critical phenomena through basic probability
  models}, volume 2025 of {\em Lecture Notes in Mathematics}.
\newblock Springer, Heidelberg, 2011.
\newblock Lecture notes from the 40th Probability Summer School held in
  Saint-Flour, 2010, \'Ecole d'\'Et\'e de Probabilit\'es de Saint-Flour.
  [Saint-Flour Probability Summer School].

\bibitem{Hug95}
B.~D. Hughes.
\newblock {\em Random walks and random environments. {V}ol. 1}.
\newblock Oxford Science Publications. The Clarendon Press, Oxford University
  Press, New York, 1995.
\newblock Random walks.

\bibitem{Hug96}
B.~D. Hughes.
\newblock {\em Random walks and random environments. {V}ol. 2}.
\newblock Oxford Science Publications. The Clarendon Press, Oxford University
  Press, New York, 1996.
\newblock Random environments.

\bibitem{Kesten75}
H.~Kesten, M.~V. Kozlov, and F.~Spitzer.
\newblock A limit law for random walk in a random environment.
\newblock {\em Compositio Math.}, 30:145--168, 1975.

\bibitem{MRS16}
A.~Matzavinos, A.~Roitershtein, and Y.~Seol.
\newblock Random walks in a sparse random environment.
\newblock {\em Electron. J. Probab.}, 21:Paper No. 72, 20, 2016.

\bibitem{Mourrat12}
J.-C. Mourrat.
\newblock On the delocalized phase of the random pinning model.
\newblock In {\em S\'eminaire de {P}robabilit\'es {XLIV}}, volume 2046 of {\em
  Lecture Notes in Math.}, pages 401--407. Springer, Heidelberg, 2012.

\bibitem{Rev75}
P.~R\'ev\'esz.
\newblock {\em Random walk in random and non-random environments}.
\newblock World Scientific Publishing Co. Pte. Ltd., Hackensack, NJ, third
  edition, 2013.

\bibitem{Seol13}
Y.~Seol.
\newblock {\em Random walks in a sparse random environment}.
\newblock ProQuest LLC, Ann Arbor, MI, 2013.
\newblock Thesis (Ph.D.)--Iowa State University.

\bibitem{Sinai82}
Y.~G. Sinai.
\newblock The limit behavior of a one-dimensional random walk in a random
  environment.
\newblock {\em Teor. Veroyatnost. i Primenen.}, 27(2):247--258, 1982.

\bibitem{Solomon75}
F.~Solomon.
\newblock Random walks in a random environment.
\newblock {\em Ann. Probability}, 3:1--31, 1975.

\bibitem{Sznitman04}
A.-S. Sznitman.
\newblock Topics in random walks in random environment.
\newblock In {\em School and {C}onference on {P}robability {T}heory}, volume
  XVII of {\em ICTP Lect. Notes}, pages 203--266. Abdus Salam Int. Cent.
  Theoret. Phys., Trieste, 2004.

\bibitem{Ton-AAP2008}
F.~L. Toninelli.
\newblock Disordered pinning models and copolymers: beyond annealed bounds.
\newblock {\em Ann. Appl. Probab.}, 18(4):1569--1587, 2008.

\bibitem{Zei2012}
O.~Zeitouni.
\newblock Random walks in random environment.
\newblock In {\em Computational complexity. {V}ols. 1--6}, pages 2564--2577.
  Springer, New York, 2012.

\end{thebibliography}


\end{document}